\documentclass[a4paper]{article}
\title{Sub-25-dimensional counterexamples to Borsuk's conjecture in the Leech lattice?}
\author{Thomas Jenrich}
\date{2025-03-13}
\begin{document}
\maketitle

\section{Abstract}

In \cite{Bor} Karol Borsuk asked whether each bounded set in the
$n$-dimensional Euclidean space can be divided into $n$+1 parts of
smaller diameter. Because it would not make sense otherwise, one usually
assumes that he just forgot to require that the whole set contains at
least two points.

The hypothesis that the answer to that question is positive became famous
under the name \emph{Borsuk's conjecture}.
Counterexamples are known for any $n\ge 64$, since 2013
(cf. \cite{Bon},\cite{Bon2},\cite{Jen},\cite{JB}).

Let $\Lambda$ be the (original, unscaled) Leech lattice, a now very well-known
infinite discrete vector set in the 24-dimensional Euclidean space.
The smallest norm of nonzero vectors in $\Lambda$ is $\sqrt{32}$. Let $M$
be the set of the 196560 vectors in $\Lambda$ having this norm.

For each $x \in M$, $-x$ is in $M$. Let $H$ be the set of all subsets of
$M$ that for each $x$ in $M$ contain either $x$ or $-x$. Each element of
$H$ has the same diameter $d = \sqrt{96}$.

For dimensions $n<24$ one can analogously construct respective $M_n$ and
$H_n$ from laminated $n$-dimensional sublattices $\Lambda_n$ of $\Lambda$.
For uniformity, let $\Lambda_{24}=\Lambda$, $M_{24} = M$ and $H_{24} = H$.
If $M_n$ is divisible into at most $n+1$ parts of diameter below $d$
then this applies to all elements of $H_n$, too. I have checked that this
is the case for all $n \le 21$. For $n$ from 22 to 24, the minimum number
of parts of diameter below $d$ that I was able to divide $M_n$ into are
25, 29 and 34, resp. The source package of this article contains a
data file encoding an element of $H_{24}=H$ that I can not divide into less than 29
parts of smaller diameter.

\section{Leech's construction of the lattice}

(See §2.31 in \cite{Lee2}) Let G be the set of the 4096 (0,1)-vectors
generated by linear combinations modulo 2 from the 12 row vectors of the
matrix C defined in §2.3 of \cite{Lee1}. A compact presentation of C:
\medskip

  100000000000011111111111

  010000000000110100011101

  001000000000101000111011

  000100000000110001110110

  000010000000100011101101

  000001000000100111011010

  000000100000101110110100

  000000010000111101101000

  000000001000111011010001

  000000000100110110100011

  000000000010101101000111

  000000000001111010001110

\medskip

$G$ is usually called the (uniqueness is up to relabeling (permutation)
of the coodinates) extended binary Golay code (named after Marcel J. E.
Golay). For any x in G, the numbers of elements of G that differ from x
in exactly 0, 8, 12, 16 or 24 coordinates are 1, 729, 2576, 729 and 1,
resp.

\medskip
The Leech lattice $\Lambda$ consists of those
$(x_1,\dots,x_{24}) \in \mathbf{Z}^{24}$ for which

$\exists a \in \{0,1\}$ and $ \exists (b_1,\dots,b_{24}) \in G $ and
$\forall i \in \{1,\dots,24\}
\ (\exists c_i \in \{0,1\} \land \exists d_i \in \mathbf{Z})$

such that
 $\sum_{i=1}^{24} c_i \equiv a\ (\rm{mod}\ 2)$ and
$ \forall i \in \{1,\dots,24\}\ x_i = a + 2b_i + 4c_i + 8d_i $.

\medskip

\section{The vectors in $\Lambda$ having squared norm 32 (the set $M$)}

The elements of $M$ are those $(x_1,\dots,x_{24})$ that can be constructed
in one of the following three ways:

(1)

For $(b_1,\dots,b_{24}) \in G$ and $p \in \{1,\dots,24\}$ :
For all $i \in \{1,\dots,24\}$
$$
x_i = \left\{\begin{array}{ll}
6b_i-3 & \mbox{if $i=p$}\\
1-2b_i & \mbox{otherwise}\end{array}\right.
$$

(2)

For $(b_1,\dots,b_{24}) \in G$ with $\sum_{i=1}^{24} b_i=8$ and
 $(c_1,\dots,c_8) \in \{0,1\}^8$ with $\sum_{i=1}^{8} c_i \equiv 0\ (\rm{mod}\ 2)$ :

Let $\{p_1,\dots,p_8\} = \{ p \in \{1,\dots,24\}\ |\ b_p=1 \}$.
For all $i \in \{1,\dots,24\}$ and all $k \in \{1,\dots,8\}$

$$
x_i = \left\{\begin{array}{ll}
2-4c_k & \mbox{if $i=p_k$}\\
0 & \mbox{if $b_i = 0$}\end{array}\right.
$$

(3)

For $p_1,p_2 \in \{1,\dots,24\}$ with $p_1<p_2$ and $s_1,s_2 \in \{0,1\}$:
For all $i \in \{1,\dots,24\}$

$$
x_i = \left\{\begin{array}{ll}
4-8s_1 & \mbox{if $i=p_1$}\\
4-8s_2 & \mbox{if $i=p_2$}\\
0 & \mbox{otherwise}\end{array}\right.
$$

\medskip
The numbers of those vectors are $4096 \cdot24=98304 $,
$759\cdot 128 = 97152 $ and ${{24}\choose{2}} \cdot 4 = 1104 $, resp.

\section{A few basic properties of $M$}

It is well-known (and can be checked computationally) that for any $x \in M$ and any $p \in \mathbf{R}$

$$
\# \{ y \in M\ |\ \langle x,y \rangle = p\} = \left\{\begin{array}{ll}
1 & \mbox{if $|p|=32$}\\
4600 & \mbox{if $|p|=16$}\\
47104 & \mbox{if $|p|=8$}\\
93150 & \mbox{if $p=0$}\\
0 & \mbox{otherwise}\end{array}\right.
$$

We are in an Euclidean space, so the squared distance of $x,y \in M$ is

$\|x-y\|^2 = \langle x-y,x-y\rangle = 2 \cdot (32-\langle x,y\rangle) $.
These cases occur:

\medskip

\begin{tabular}{ |c|r|r|r|r|r|r|r|r|r|r|r|r|r|r|r|r|r|r|r| }
\hline
 $\langle x,y\rangle$&32&16&8&0&-8&-16&-32\\
\hline
 $\|x-y\|^2$&0&32&48&64&80&96&128\\

\hline
\end{tabular}

\medskip

Thus, the diameter of $M$ is $\sqrt{128}$ and for any element of $H$
the diameter is $d = \sqrt{96}$ and the (easier to calculate) smallest
inner product is -16.

\section{On other information sources}

Because of its remarkable properties there is a lot of literature on the
Leech lattice. A very early article that described  many of
the symmetries of the lattice (and of $M$, resp.) is \cite{Con}.

The book \cite{CS2} is often recommended as standard source on this
topic. (Personally, I only have read the article \cite{CS1} which
seemingly became a part of \cite{CS2} later).

Be aware that it is not unusual to scale the Leech lattice by a factor,
e.g. $\sqrt{1/8}$ or $\sqrt{1/32}$, in order to get the minimal positive
norm 2 or 1, resp., which let it fit into certain frameworks.

\section{Parameters of the Leech lattice and its laminated sublattices}

The conditions in the table below have been derived from Figure 4 in
\cite{CS1}. As in that figure, the other laminated sublattices in
dimensions 11, 12 and 13 have been omitted.

\medskip

\begin{tabular}{ | r | l | l | l | r | r | c | c | }
\hline
n & Latt. & $\cong$ & added condition  & $\# M_n$ & $(\pm 4,0)$ & $(\pm 2,0)$ & $(\pm 3,\pm 1)$\\
\hline
24& $\Lambda_{24}$ & & & 196560 &1104&$97152=759\cdot 128$&$98304=24 \cdot 4096$\\
23& $\Lambda_{23}$ & & $x_{24}=x_{23}$&93150&926&$47168=737 \cdot 64 $&$45056=22 \cdot 2048$\\
22& $\Lambda_{22}$ & & $x_{23}=x_{22}$ &  49896 &840&$27552=861 \cdot 32$&$21504=21 \cdot 1024$\\
21& $\Lambda_{21}$ & & $x_{22}=0$ &  27720 &840&$26880=210 \cdot 128$&0\\
20& $\Lambda_{20}$ & & $x_{21}=0$ &  17400 &760&$16640=130 \cdot 128$&0\\
19& $\Lambda_{19}$ & & $x_{20}=0$ &  10668 &684&$9984=78 \cdot 128$&0\\
18& $\Lambda_{18}$ & & $x_{19}+x_{18}+x_{17}=0$ & 7398 &486&$6912=54 \cdot 128$&0\\
17& $\Lambda_{17}$ & & $x_{19}=0$  &  5346 &482&$4864=38 \cdot 128$&0\\
16& $\Lambda_{16}$ & BW16 & $x_{18}=0$ &   4320 &480&$3840=30 \cdot 128$&0\\
15& $\Lambda_{15}$ & & $x_{16}=0$ &   2340 &420&$1920=15 \cdot 128$&0\\
14& $\Lambda_{14}$ & & $x_{15}+x_{14}+x_{13}=0$ & 1422 &270&$1152=9 \cdot 128$&0\\
13& $\Lambda_{13}^{max}$ & & $x_{15}=0$ & 906 &266&$640=5 \cdot 128$&0\\
12& $\Lambda_{12}^{max}$ & & $x_{14}=0$ & 648 &264&$384=3 \cdot 128$&0\\
11& $\Lambda_{11}^{max}$ & & $x_{12}=x_{11}$  & 438 &182&$256=2 \cdot 128$&0\\
10& $\Lambda_{10}$ & & $x_{11}=x_{10}$  & 336 &144&$192=3 \cdot 64$&0\\
9& $\Lambda_9 $ & & $x_{10}=0$ & 272 &144&128&0\\
8& $\Lambda_8 $ & $E_8$ & $x_{9}=0$ & 240 &112&128&0\\
7& $\Lambda_7 $ & $E_7$ & $x_8=x_7$  & 126 &62&64&0\\
6& $\Lambda_6 $ & $E_6$ & $x_7=x_6$  & 72 &40&32&0\\
5& $\Lambda_5 $ & $D_5$ & $x_{6}=0$ & 40 &40&0&0\\
4& $\Lambda_4 $ & $D_4$ & $x_{5}=0$ & 24 &24&0&0\\
3& $\Lambda_3 $ & $D_3$,$A_3$ & $x_{4}=0$ & 12 &12&0&0\\
2& $\Lambda_2 $ & $A_2$ & $ x_3+x_2+x_1=0 $ & 6 &6&0&0\\
1& $\Lambda_1 $ & $A_1$,$\mathbf{Z}$ & $ x_3=0 $ & 2 &2&0&0\\

\hline
\end{tabular}

\medskip

Inner products in $M_n$ are never $8$ or $-8$ for $n \le 8$, are always
$32$ or $-32$ for $n=1$.

\section{Dividing $M_n$ and some elements of $H_n$ - efforts and results}

Because of the strong regularity of the $M_n$, a rather mathematical
treatment seems to be not just preferable but also applicable with success.
But I gave it up after not reaching (all of) the results from actual
computations.
So, the results stated below have been found and independently re-checked
(almost) completely computationally.

For any subset $X$ of $M$ let $\Gamma_X$ be the graph whose vertices are the
vectors in $X$, adjacent if and only if their inner product is not
above $-16$ . A proper (i.e., conflict-free) coloring of $\Gamma_X$ with $c$ colors
corresponds to a division of $X$ into $c$ parts of diameter below $d$.

It is well known that many problems in mathematics and computer science
(of theoretical as well as of practical nature) can be, have been and
are translated into tasks to find proper colorings of the vertices of
graphs.
There are many papers describing and comparing graph coloring algorithms
and the results of their application to certain concrete graphs, mostly
from benchmark collections published in the second implementation
challenge by the DIMACS (Center for Discrete Mathematics and Theoretical
Computer Science) in 1992-1993 . Unfortunately, the described
implementations are not (freely) downloadable or not applicable to the
bigger of the graphs considered in this article (on a computer with
4 GB RAM). Generally, suitable implementations seem to be hard to get.
As you can read in older versions of this article I even took the
detour of converting the problem into input for SAT solvers. And I also
tried to employ programs performing Mixed Integer Linear Programming. Later I found a bundle of graph coloring programs
implementing different algorithms, exact and heuristic ones. As it turned
out, the program implementing the heuristic algorithm known as TABUCOL
(1987 described by A. Hertz and D. de Werra) delivered the best results
in an acceptable time.
After some modifications of the source code I could use it to treat $M_n$
for any $n \le {23}$ and to treat any set in each $H_n$ on my PC with 4
GB RAM. An important part of the TABUCOL algorithm is randomization. The
implementing program uses numbers from a pseudo-randomizing generator
depending on a seed value given within the command line. So I did
multiple runs with different seed values in order to improve already
found results.
\medskip

In order to be able to handle the bigger of the considered graphs
(at all or more efficiently), I have found in and removed from the
given graph some large independent sets before applying that program.
To estimate (sets of) such independent sets, the following thoughts and
insights have been helpful (without guaranteeing not to miss a
perhaps possible division of the complete set into fewer parts):

Consider $x, y \in M$ with $\langle x,y \rangle= (-8)$. The sum of $x$
and $y$ is in $\Lambda$ but has the norm $\sqrt{48}$ (and is therefore
not in $M$). The smallest inner product of vectors in
$\{ z \in M : \langle x+y, z \rangle \ge 16 \}$ is $-8$. The size of
this set is 11730 which seems to be the maximum size of subsets of $M$
with smallest inner product above $-16$.

Though $25 \times 11730 > \#M$, it is not possible to cover $M$ by 25 of
such sets because of their ``roundness''. Wherever the centers of these
sets are placed, some regions remain uncovered.
In particular, it is not possible to have 25 such centers with the
same pairwise inner product (as like as the corners of a 24-dimensional
regular simplex), which for vectors of their norm would have to be $-2$.
The closest actually appearing inner products are 0 and -8.

For dimension $n$ from 9 to 23, one can analogously consider
$x,y \in M_n$ with $\langle x,y \rangle= (-8)$ and the set
 $\{ z \in M_n : \langle x+y, z \rangle \ge 16 \}$. Clearly, such sets
are smaller in lower dimensions than in higher ones. And the size can
depend on $x+y$ because not all $M_n$ are that symmetric as $M$ is. But I
think that the largest subsets of $M_n$ with smallest inner product
above $-16$ can be constructed this way.

\subsection{Divisions of the $M_n$ (summary)}

This table contains for dimensions $n$ from 1 to 24 the smallest number
of parts of diameter below $d$ that I have been able to divide $M_n$ into:

\medskip

\begin{tabular}{ |r|r|r|r|r|r|r|r|r|r|r|r|r|r|r|r|r|r|r| }
\hline
 1&2&3&4&5&6&7&8, ..., 12&13&14&15&16&17&18, 19&20&21&22&23&24\\
\hline
 2&3&4&5&6&7&8&9&10&11&13&15&16&17&20&22&25&29&34\\
\hline
\end{tabular}

\medskip

Compared to the first two versions of this article, there are reductions
for dimensions 14, 16, 17, 18, 19, 20, 22, 23, 24. Concerning Borsuk's
conjecture, 22 is the smallest remaining interesting dimension.

\subsection{A potential 24-dimensional counterexample to Borsuk's conjecture}

The elements of $H_{22}$ and $H_{23}$ that I have investigated can be
divided into 22 and 24, resp., parts of diameter below $d$, are thus no
counterexamples to Borsuk's conjecture.

Of course, this does not mean that this applies to all elements of $H_{22}$
and $H_{23}$.

Within $H_{24}$, on the other hand, I have found a set that I tried but
failed to divide into less than 29 parts of diameter below $d$.
Thus, I consider that set as a potential but currently unverified
counterexample to Borsuk's conjecture.

In order to make that set available for investigation by others,
the source package of this article contains the binary data file
H24S1.DAT that encodes a sequence of the 98280 vectors in that set,
in the order and coding given in this table:

\medskip

\begin{tabular}{ |r|l|r|r|r|r| }
\hline
 number & type & 0 & 1 & 2 & 3 \\
\hline
 552&$(\pm 4, 0)$& -4 & 0 & 4 & \\
\hline
 48576&$(\pm 2, 0)$& -2 & 0 & 2 &\\
\hline
 49152&$(\pm 3, \pm 1)$ & -3 & -1 & 1 & 3 \\
\hline
\end{tabular}

\medskip

The result $(y_1, ..., y_{24})$ of the component-wise encoding of
a vector is represented by 6 bytes, where the $n$-th byte equals
$ 64 \cdot y_{4n}  + 16 \cdot y_{4n-1} + 4 \cdot y_{4n-2} + y_{4n-3}$.

\vspace{0.1in}

Author's eMail address: thomas.jenrich@gmx.de

\end{document}